\documentclass[12pt]{article}
\usepackage{amsfonts}
\usepackage{amssymb}

\input{tcilatex}
\usepackage{polski}
\begin{document}

\begin{center}
\textbf{Cauchy's functional equation and extensions: Goldie's equation and
inequality, the Go\l \k{a}b-Schinzel equation and Beurling's equation}\\[0pt]
\textbf{by\\[0pt]
N. H. Bingham and A. J. Ostaszewski}

\bigskip

\textit{To Ranko Bojani\'{c} on his }$90^{th}$ \textit{birthday.}
\end{center}

\bigskip

\noindent \textbf{Abstract. }The Cauchy functional equation is not only the
most important single functional equation, it is also central to regular
variation. Classical Karamata regular variation involves a functional
equation and inequality due to Goldie; we study this, and its counterpart in
Beurling regular variation, together with the related Go\l \k{a}b-Schinzel
equation.


\bigskip

\section{Introduction}

We are concerned with the most important single functional equation, the
\textit{Cauchy Functional Equation $(CFE)$}

\begin{equation}
K(x+y)=K(x)+K(y),\text{\qquad }k(xy)=k(x)k(y),  \tag{$CFE$}
\end{equation}%
to give both the additive and multiplicative versions. For background, see
the standard work by Kuczma [Kuc]. This is known to be crucial to the theory
of regular variation, in both its Karamata form (see Ch. 1 of [BinGT], BGT
below) and its Bojani\'{c}-Karamata/de Haan form (BGT Ch. 3, [BojK]). A
close study of these involves a certain functional equation ([BinG], BGT),
which we call here the \textit{Goldie} functional equation ($(GFE)$ -- see
below) \footnote{%
The equation occurs first in joint work by the first author and Goldie; the
first author is happy to confirm that the argument is in fact due to Goldie,
whence the name.}. One of the themes of Kuczma's book is the interplay
between functional equations and inequalities; he focusses particularly on
the Cauchy functional equation and Jensen's inequality. Even more closely
linked to $(CFE)$ is the functional inequality of \textit{subadditivity},
and this has its counterpart in \textit{Goldie's functional inequality} ($%
(GFI)$ -- see below).

Closely related to the Karamata theory of slow and regular variation is the
theory of \textit{Beurling slow and regular variation}. This has an odd
history. Beurling (who was a perfectionist) never published it (it is not
mentioned in the two volumes of his Collected Works). He introduced it in
lectures in the 1940s for use in his \textit{Tauberian theorem}. Beurling's
Tauberian theorem appeared in the 1972 papers of Moh [Moh] and Peterson
[Pet]; it was used by the first author [Bin] in 1981 in probability theory.
Beurling slow variation has been in use ever since (see BGT \S 2.11 and
[Kor] for background and references). Beurling regular variation was
introduced in our recent paper [BinO3]. Here it emerged that the Beurling
theory of slow and regular variation \textit{includes} the Karamata theory
(despite Beurling slow variation having previously appeared to be a minor
topic within the Karamata theory). The role of $(CFE)$ and $(GFE)$ is played
in this Beurling context by a functional equation, which we call here the
\textit{Beurling functional equation} ($(BFE)$ -- see below). This
functional equation is a special case of one independently introduced by Acz%
\'{e}l and by Go\l \k{a}b and Schinzel [GolS] in 1959, and (without
knowledge of Beurling's work) was called the \textit{Go\l \k{a}b-Schinzel
functional equation} in Acz\'{e}l and Dhombres [AczD]; see \S 6.2 for more
current information on this literature and associated applications. It has
recently been studied by the second author [Ost2] for its relation with
uniform convergence in the context of Beurling regular variation; here we
unify these two lines of work.

The theme of the present paper is that one begins with the functional
inequality, imposes a suitable side-condition (which serves to `give the
inequality the other way') and deduces the corresponding functional
equation, which under suitable conditions one is able to solve. The
functional equation and functional inequality we have in mind are those
mentioned above:
\begin{equation}
F^{\ast }(u+v)\leq e^{\rho v}F^{\ast }(u)+F^{\ast }(v)\text{ \qquad }%
(\forall u,v\in \mathbb{R})  \tag{$GFI$}
\end{equation}%
(BGT\ (3.2.5), cf. (3.0.11) -- for $F^{\ast }$ and $F$ see below), and
\begin{equation}
F^{\ast }(u+v)=e^{\rho v}K(u)+F^{\ast }(v)\text{ \qquad }(\forall u,v\in
\mathbb{R})  \tag{$GFE$}
\end{equation}%
(see BGT\ (3.2.7) -- for the relationship here for $\rho \neq 0$ between $%
F^{\ast }$ and its \textit{Goldie kernel} $K$, indeed in greater generality,
see Th. 3 in \S 2 below; cf. [BinO4, Prop. 1] for the \textit{additive kernel%
} in the case $\rho =0,$ which reduces $(GFI)$ to subadditivity on $\mathbb{R%
}_{+}$, the context there). $(GFI)$ captures an asymptotic relation in
functional form, and so is key to establishing the Characterization Theorem
of regular variation (BGT \S 1.4). Our focus here is on the extent to which
the universal quantifiers occurring in the functional inequalities and
functional equations under study can be weakened, in the presence of
suitable side-conditions. The prototypical side-condition here is the
\textit{Heiberg-Seneta condition }%
\begin{equation}
\lim \sup_{u\downarrow 0}F(u)\leq 0,  \tag{$HS(F)$}
\end{equation}%
due to Heiberg [Hei] in 1971, and Seneta [Sen] in 1976 (BGT, Th. 3.2.4).
This condition, which is best possible here, is what is needed to reduce $%
(GFI)$ to $(GFE)$.

Two related matters occur here. One is the question of \textit{quantifier
weakening} above. This, together with $(HS),$ hinges on the algebraic nature
of the set on which one can assert equality. The second, \textit{automatic
continuity}, relates to the extent to which a solution of $(GFE)$ is
continuous (and hence easily of standard form -- see BGT Ch. 3), or (in the
most important case $\rho =0$) an \textit{additive} function becomes \textit{%
continuous}, and so \textit{linear}. This is the instance of the important
subject of \textit{automatic continuity }relevant here. Automatic continuity
has a vast literature, particularly concerning homomorphisms of Banach
algebras, for which see [Dal1] and [Dal2]. See also Helson [Hel] for Gelfand
Theory, Ng and Warner [NgW] and Hoffman-J\o rgensen [HJ]. The crux here is
the \textit{dichotomy} between additive functions with a hint of regularity,
which are then linear, and those without, which are pathological. For
background and references on dichotomies of this nature, Hamel pathologies
and the like, see [BinO2].

One of our themes here and in [BinO4] is \textit{quantifier weakening}: one
weakens a universal quantifier $\forall $ by thinning the set over which it
ranges. In what follows we will often have two quantifiers in play, and will
replace \textquotedblleft $\forall u\in \mathbb{A}$\textquotedblright\ by
\textquotedblleft $(u\in \mathbb{A)}$\textquotedblright , etc. This
convention, convenient here, is borrowed from mathematical logic.

One theme that this paper and [BinO4] have in common is the great debt that
the subject of regular variation, as it has developed since [BinG] and BGT,
owes to the Goldie argument. It is a pleasure to emphasize this here. This
argument originated in a study of \textit{Frullani integrals}, important in
many areas of analysis and probability ([BinG I, II.6]; cf. BGT \S 1.6.4,
[BinO4, \S 1]).

\section{Generalized Goldie equation}

We begin by generalizing $(GFE)$ by replacing the exponential function on
the right by a more general function $g,$ the \textit{auxiliary function. }%
We further generalize by \textit{weakening the quantifiers}, allowing them
to range over a set $\mathbb{A}$ smaller than $\mathbb{R}$. It is
appropriate to take $\mathbb{A}$ as a dense (additive) subgroup. The
functional equation in the result below, written there $(G_{\mathbb{A}})$,
may be thought of as the second form of the Goldie functional equation
above. As we see in Theorem 1 below, the two coincide in the principal case
of interest -- compare the insightful Footnote 3 of [BojK]. The notation $%
H_{\rho }$ below is from BGT \S 3.1.7 and 3.2.1 implying $H_{0}(t)\equiv t.$
The identity $uv-u-v+1=(1-u)(1-v)$ gives that $(1-e^{-\rho x})/\rho $ is
\textit{subadditive} on $\mathbb{R}_{+}:=(0,\infty )$ for $\rho \geq 0,$ and
\textit{superadditive} on $\mathbb{R}_{+}$ for $\rho \leq 0.$

\bigskip

\noindent \textbf{Theorem 1.} \textit{For }$g$ \textit{with }$g(0)=1,$
\textit{if }$K\neq 0$\textit{\ satisfies }$(G_{\mathbb{A}})$\textit{\ below
with }$\mathbb{A}$\ \textit{a dense subgroup:}%
\begin{equation}
K(u+v)=g(v)K(u)+K(v),\text{\qquad }(u,v\in \mathbb{A})
\tag{$G_{\QTR{Bbb}{A}}$}
\end{equation}%
-- \textit{then }%
\[
\mathbb{A}_{g}:=\{u\in \mathbb{A}:\mathit{\ }g(u)=1\}
\]%
\textit{\ is an additive subgroup on which }$K$\textit{\ is additive, and
for some constant }$\kappa $%
\begin{equation}
K(t)\equiv \kappa (g(t)-1)\text{ for }t\in \{0\}\cup \mathbb{A}\text{%
\TEXTsymbol{\backslash}}\mathbb{A}_{g}.  \tag{*}
\end{equation}%
\textit{\ For }$\mathbb{A=R}$\textit{\ and }$g$\textit{\ locally bounded at }%
$0$\textit{\ with }$g\neq 1$\textit{\ except at }$0:$\textit{\ }$g(x)\equiv
e^{-\rho x}$ \textit{for some constant }$\rho \neq 0,$ \textit{and so }$%
K(t)\equiv \kappa \rho H_{\rho }(t),$ \textit{where}%
\[
H_{\rho }(t):=(1-e^{-\rho t})/\rho .
\]

\bigskip

\noindent \textit{Proof.} Recall that
\[
\mathcal{N}_{K}:=\{x\in \mathbb{A}:K(x+a)=K(x)+K(a)\text{ }(\forall a\in
\mathbb{A})\}
\]%
is the \textit{Cauchy nucleus} of $K$ -- see [Kuc, \S 18.5], and is either
empty or a subgroup (for a proof see [Kuc, Lemma 18.5.1], or the related
[AczD, Ch. 6, proof of Th. 1). If $x\in \mathcal{N}_{K},$ choosing $a\in
\mathbb{A}$ with $K(a)\neq 0$ yields $g(x)=1$ from%
\[
K(a+x)=K(a)+K(x)=g(x)K(a)+K(x).
\]%
Conversely, for $v\in \mathbb{A}_{g}$ and any $v\in \mathbb{A}$%
\[
K(u+v)=K(u)+K(v)\text{ }
\]%
so $v\in \mathcal{N}_{K}$: $\mathbb{A}_{g}=\mathcal{N}_{K}.$ By assumption $%
0\in \mathbb{A}_{g},$ so in particular $K\ $is additive on $\mathbb{A}_{g},$
and $K(0)=0$.

Continue now as in BGT Lemma 3.2.1: for $u,v\in \mathbb{A}$\TEXTsymbol{%
\backslash}$\mathbb{A}_{g}$ distinct:%
\[
g(v)K(u)+K(v)=K(u+v)=K(v+u)=g(u)K(v)+K(u).
\]%
So%
\[
K(u)[g(v)-1]=K(v)[g(u)-1]:\qquad \frac{K(u)}{g(u)-1}=\frac{K(v)}{g(v)-1}%
=\kappa ,
\]%
say; so, for some constant $\kappa ,$
\[
K(u)=\kappa \lbrack g(u)-1]
\]%
on $\mathbb{A}\backslash \mathbb{A}_{g}$, proving ($\ast )$ in this case.
Although we assumed $u\neq 0,$ we still have $0=K(0)=\kappa \lbrack g(0)-1],$
completing the proof of ($\ast $).

Substitution (for $u,v\in \mathbb{A}\backslash \mathbb{A}_{g}$ provided $%
u+v\in \mathbb{A}\backslash \mathbb{A}_{g}$) yields first%
\[
\kappa \lbrack g(u+v)-1]=\kappa g(v)[g(u)-1]+\kappa \lbrack g(v)-1],
\]%
and then, for $\kappa \neq 0,$ the Cauchy exponential equation%
\begin{equation}
g(u+v)=g(v)g(u),  \tag{$CEE$}
\end{equation}%
if $\mathbb{A}_{g}=\{0\};$ so if $\mathbb{A=R}$, local boundedness yields $%
g(x)\equiv e^{-\rho x}$ for some $\rho $ (see [AczD, Ch. 3], or [Kuc, \S %
13.1]). If $\kappa =0$ above, then $K(x)\equiv 0.$ $\square $

\bigskip

\noindent \textit{Remarks. 1. }In Theorem 1 above the additive reals act on
the domain of the unknown function $K.$ Generalizations are possible to
other group actions and will rely on the auxiliary function $g$ being a
group homomorphism as in $(CEE)$ above. There is more to be said here; we
hope to this elsewhere.

\noindent \textit{2.} Recall that $f$ satisfies the \textit{Mikusi\'{n}ski
equation }if%
\begin{equation}
f(x+y)=f(x)+f(y)\qquad \text{if }f(x+y)\neq 0;  \tag{$Mik$}
\end{equation}%
such a function is necessarily additive, for which see [AczD, Ch. 6 Th. 1].
The argument above identifies $g$ from $(CEE)$ when $u,v\in \mathbb{A}%
\backslash \mathbb{A}_{g}$, provided $u+v\in \mathbb{A}\backslash \mathbb{A}%
_{g};$ so for $g>0,$ the condition $g(u+v)\neq 1$ is equivalent to $\log
g\neq 0,$ which means that $\log g$ satisfies $(Mik)$ and so $g(x)=\exp f(x)$
for some additive (possibly pathological) function $f$.

\noindent \textit{3. }Above, for $g$ Baire/measurable, by the Steinhaus
subgroup theorem (see e.g. [BinO2, Th. S] for its general combinatorial
form), $\mathbb{A}_{g}=\mathbb{R}$ iff $\mathbb{A}_{g}$ is non-negligible,
in which case $K\ $is additive. The additive case is studied in [BinO4] and
here we have passed to $\mathbb{A}_{g}=\{0\}$ as a convenient context. More,
however, is true. As an alternative to the last remark, for $\mathbb{A}_{g}$
negligible: by the Fubini/Kuratowski-Ulam Theorem (for which see [Oxt, Ch.
14-15]), the equation $(CEE)$ above holds for quasi almost all $(u,v)\in
\mathbb{R}^{2}$; consequently, by a theorem of Ger (see [Ger], or [Kuc,
Th.18.71]), there is a homomorphism on $\mathbb{A}$ `essentially extending' $%
\log g$ to $\mathbb{A}.$ From here, again for $g$ Baire/measurable, $%
g(x)=e^{-\rho x}$ for some $\rho .$

\bigskip

In Theorem 2 below, there is no quantifier weakening to $\mathbb{A}$ and so
we need $(G_{\mathbb{R}})$ in place of $(G_{\mathbb{A}}).$ It will be
convenient in what follows to write `positive' for functions to mean
`positive on $\mathbb{R}_{+}$', unless otherwise stated.

\bigskip

\noindent \textbf{Theorem 2.} \textit{If both }$K$ \textit{and }$g$ \textit{%
in }$(G_{\mathbb{R}})$\textit{\ are positive with }$g\neq 1$ \textit{except
at }$0$\textit{,\ then either }$K\equiv 0,$ \textit{or both are continuous,
and }$g(x)\equiv e^{-\rho x}$ \textit{for some }$\rho \neq 0.$

\bigskip

\noindent \textit{Proof}\textbf{.} Writing $w=u+v,$ one has%
\[
K(w)-K(v)=g(v)K(w-v),
\]%
so $K$ is strictly increasing and so continuous at some point $y>0$ say. But
for any $h$%
\[
K(y+h)-K(y)=g(y)K(h),
\]%
and so, since $g(y)>0$, $K$ is continuous at $0,$ as $K(0)=0.$ Hence $K$ is
continuous at any point $t>0,$ since $g(t)>0$ and%
\[
K(t+h)-K(t)=g(t)K(h).
\]%
Take any $u>0.$ Fix $w>u$, so that $K(w-u)>0.$ Then, since%
\[
g(u)=[K(w)-K(u)]/K(w-u),
\]%
and the right-hand side is continuous in $u$ for $u>0,$ the function $g$ is
continuous for $u>0.$ Finally, as in Th. 1, $K(x)=\kappa \lbrack g(x)-1]$
for all $x,$ for some constant $\kappa ;$ if $\kappa \neq 0,$ then $g$
satisfies $(CEE)$ and is continuous on $\mathbb{R}_{+},$ so again the
function $g$ is $e^{-\rho x}$ and $K$ is continuous. $\square $

\bigskip

\noindent \textit{Remark. }For $g\equiv 1$ in $(G_{\mathbb{R}}),$ the proof
of Theorem 2 shows that a positive additive function is continuous -- a weak
form of Darboux's Theorem on the continuity of bounded additive functions.

\bigskip

In\textbf{\ }$(G_{\mathbb{R}})$ above for $x,\rho \geq 0$ one has $%
g(x)=e^{-\rho x}\leq 1$ on $\mathbb{R}_{+}$; generally, if $g(x)\leq 1$ on $%
\mathbb{R}_{+}$ and $K\ $positive satisfies $(G_{\mathbb{R}}),$ then for $%
u,v\geq 0$%
\[
K(u+v)\leq K(u)+K(v),
\]%
and so $K$ is subadditive on $\mathbb{R}_{+}$.

\bigskip

We now prove a converse -- our main result. Here, in the context of
subadditivity, the important role of the \textit{Heiberg-Seneta condition},
discussed in \S 1, is performed by a weaker side-condition: \textit{%
right-continuity at} $0$, a consequence, established in [BinG] -- see also
BGT \S 3.2.1 and [BinO3]. A further quantifier weakening occurs in (ii)
below.

\bigskip

\noindent \textbf{Theorem 3 (Generalized Goldie Theorem).} \textit{If for }$%
\mathbb{A}$\textit{\ a dense subgroup,}

\noindent (i)\textit{\ }$F^{\ast }:$\textit{\ }$\mathbb{R\rightarrow R}$
\textit{is positive and subadditive with }$F^{\ast }(0+)=0;$

\noindent (ii) $F^{\ast }$ \textit{satisfies the weakened Goldie equation}%
\[
F^{\ast }(u+v)=g(v)K(u)+F^{\ast }(v)\text{\qquad }(u\in \mathbb{A})(v\in
\mathbb{R}_{+})
\]%
\textit{for some non-zero }$K$\textit{\ satisfying }$(G_{\mathbb{A}})$%
\textit{\ with }$g$ \textit{continuous and }$\mathbb{A}_{g}=\{0\}$\textit{;}

\noindent (iii)\textit{\ }$F^{\ast }$ \textit{extends }$K$ \textit{on }$%
\mathbb{A}$:\textit{\ }%
\[
F^{\ast }(x)=K(x)\text{ for }x\in \mathbb{A},
\]%
\textit{so that in particular }$F^{\ast }$\textit{\ satisfies }$(G_{\mathbb{A%
}}),$ \textit{and indeed}%
\[
F^{\ast }(u+v)=g(v)F^{\ast }(u)+F^{\ast }(v)\text{ for }u\in \mathbb{A},v\in
\mathbb{R}_{+};
\]%
-- \textit{then for some }$c>0,\rho \geq 0$%
\[
g(x)\equiv e^{-\rho x}\text{ and }F^{\ast }(x)\equiv cH_{\rho
}(x)=c(1-e^{-\rho x})/\rho .
\]

\bigskip

\noindent \textit{Proof. }Put%
\[
\gamma (x)=\int_{0}^{x}g(t)dt:\text{\qquad }\gamma ^{\prime }(x)=g(x).
\]%
By continuity of $g$ and Th. 1, $K(u+)=K(u)$ for all $u\in \mathbb{A}$, and
so $K(0+)=0.$ Furthermore, note that $F^{\ast }$ is right-continuous on $%
\mathbb{A}$ (and $F^{\ast }(u+)=K(u)$ on $\mathbb{A}$), and on $\mathbb{R}%
_{+}$%
\[
\lim \sup_{v\downarrow 0}F^{\ast }(u+v)\leq F^{\ast }(u)+F^{\ast
}(0+)=F^{\ast }(u).
\]

Now proceed as in the Goldie proof -- see e.g. BGT \S 3.2.1. (This uses the
sequence $s_{n}=n\delta ,$ rather than the Beck sequence of \S 3 below which
is not appropriate here, but see below in Theorem 7 for a Beck-sequence
adaptation of the current argument.) For any $u,u_{0}$ with $u_{0}\in
\mathbb{A}$ and $u_{0}>0$, define $i=i(\delta )\in \mathbb{Z}$ for $\delta
>0 $ so that $\left( i-1)\right) \delta \leq u<i\delta ,$ and likewise for $%
u_{0}$ define $i_{0}(\delta ).$ Also put
\[
c_{0}:=K(u_{0})/[g(u_{0})-1].
\]%
As $m\delta \in \mathbb{A}$,%
\[
F^{\ast }(m\delta )-F^{\ast }((m-1)\delta )=g((m-1)\delta )K(\delta ),
\]%
so that on summing%
\begin{equation}
F^{\ast }(i(\delta )\delta )=K(\delta )\sum_{m=1}^{i}g((m-1)\delta ),
\tag{**}
\end{equation}%
as $F^{\ast }(0)=0.$ Note that as $\delta \rightarrow 0,$%
\begin{equation}
\delta \sum_{m=1}^{i}g((m-1)\delta )\rightarrow \int_{0}^{u}g(x)dx
\tag{$RI$}
\end{equation}%
(for `Riemann Integral'). Assuming we may take limits as $\delta \rightarrow
0$ through positive $\delta \in \mathbb{A}$ with $K(\delta )\neq 0,$ we then
have%
\[
\frac{F^{\ast }(i(\delta )\delta )}{F^{\ast }(i_{0}(\delta )\delta )}=\frac{%
K(\delta )}{K(\delta )}\frac{\sum_{m=1}^{i}g((m-1)\delta )}{%
\sum_{m=1}^{i_{0}}g((m-1)\delta )}=\frac{\delta \sum_{m=1}^{i}g((m-1)\delta )%
}{\delta \sum_{m=1}^{i_{0}}g((m-1)\delta )}\rightarrow \frac{\gamma (u)}{%
\gamma (u_{0})}.
\]%
So by right-continuity at $u_{0},$%
\[
\lim F^{\ast }(i_{0}(\delta )\delta )=F^{\ast
}(u_{0})=K(u_{0})=c_{0}[g(u_{0})-1].
\]%
So%
\[
F^{\ast }(i(\delta )\delta )\rightarrow \gamma (u)\cdot F^{\ast
}(u_{0})/\gamma (u_{0}).
\]%
As before, as $u_{0}\in \mathbb{A},$%
\begin{eqnarray*}
F^{\ast }(u) &\geq &\lim \sup F^{\ast }(i(\delta )\delta )=\gamma (u)\cdot
F^{\ast }(u_{0})/\gamma (u_{0}) \\
&=&\gamma (u)K(u_{0})/\gamma (u_{0})=\gamma (u)c_{0}[g(u_{0})-1]/\gamma
(u_{0}).
\end{eqnarray*}%
Put%
\[
c_{1}=c_{0}[g(u_{0})-1]/\gamma (u_{0}).
\]%
Now specialize to $u\in \mathbb{A},$ on which, by above, $F^{\ast }$ is
right-continuous. Letting $i(\delta )\delta \in \mathbb{A}$ decrease to $u,$
the inequality above becomes an equation:%
\[
K(u)=F^{\ast }(u)=\gamma (u)c_{0}[g(u_{0})-1]/\gamma (u_{0})=c_{1}\gamma
(u).
\]%
This result remains valid with $c_{1}=0$ if $K(\delta )=0$ for $\delta \in
\mathbb{A}\cap I$ for some interval $I=(0,\varepsilon ),$ as then $F^{\ast
}(u)=0$ by right-continuity, because $F^{\ast }(i(\delta )\delta )=0$ for $%
\delta \in \mathbb{A}\cap I,$ by (**).

Now for \textit{arbitrary} $u,$ taking $v\uparrow u$ with $v\in \mathbb{A},$
we have (as $u-v>0)$ that%
\begin{eqnarray*}
F^{\ast }(u) &=&F^{\ast }(u-v+v)=K(v)g(u-v)+F^{\ast }(u-v)\qquad \text{(by
(ii), as }v\in \mathbb{A}\text{)} \\
&=&c_{1}\gamma (v)g(u-v)+F^{\ast }(u-v)\rightarrow c_{1}\gamma (u),
\end{eqnarray*}%
by continuity of $\gamma .$ Thus for \textit{all} $u,$%
\[
F^{\ast }(u)=c_{1}\gamma (u).
\]%
Thus by Theorem 1, for some $\kappa $%
\[
c_{1}\gamma (u)=F^{\ast }(u)=K(u)=\kappa \lbrack g(u)-1]
\]%
on $\mathbb{A}$. So, by density and continuity on $\mathbb{R}_{+}$,
\[
\kappa \lbrack g(u)-1]=c_{1}\gamma (u).
\]%
So $g$ is indeed differentiable; differentiation now yields%
\[
c_{1}g(u)=\kappa g^{\prime }(u).
\]%
If $\kappa =0,$ then $K(u)\equiv 0,$ contrary to assumptions. So%
\[
g^{\prime }(u)=(c_{1}/\kappa )g(u),
\]%
and so with $\rho :=-c_{1}/\kappa $%
\[
g(u)=e^{-\rho u}\text{ and }\gamma (u)=H_{\rho }(u):\text{\qquad }F^{\ast
}(u)=c_{1}\gamma (u)=c_{1}[1-e^{-\rho u}]/\rho .
\]%
Finally, as $(1-e^{-\rho x})/\rho $ is subadditive iff $\rho \geq 0$ (cf.
before Theorem 1), $c_{1}>0.$ $\square $

\bigskip

\noindent \textit{Remark. }From the passage from the Riemann sum in $(RI)$
to the Riemann integral, we see the origin of the otherwise surprising
feature: that we obtain automatic differentiability from continuity in
several of the arguments below.

\bigskip

\noindent \textbf{Theorem 4.} \textit{If }$g,K\ $\textit{are positive, }$%
F^{\ast }$\textit{\ is subadditive on }$\mathbb{R}_{+}$\textit{\ with }$%
F^{\ast }(0+)=0,$\textit{\ and}
\[
F^{\ast }(u+v)=g(v)K(u)+F^{\ast }(v)\text{\qquad }(u\in \mathbb{A})(v\in
\mathbb{R}_{+})
\]%
-- \textit{then }$F^{\ast }$\textit{\ is increasing and continuous on }$%
\mathbb{R}_{+}$, \textit{and so }$g$\textit{\ is continuous on }$\mathbb{R}%
_{+}.$ \textit{In particular, the continuity assumed in Theorem 3 above is
implied by positivity of both }$g$ \textit{and }$K$\textit{.}

\bigskip

\noindent \textit{Proof. }Since
\[
F^{\ast }(v+u)-F^{\ast }(v)=g(v)K(u)\text{\qquad }(u\in \mathbb{A})(v\in
\mathbb{R}_{+}),
\]%
then for $u>0$ and $u\in \mathbb{A}$
\[
F^{\ast }(v+u)>F^{\ast }(v)\text{\qquad }(v\in \mathbb{R}).
\]%
So letting $u\downarrow 0$ through $\mathbb{A}$,%
\[
F^{\ast }(v)\leq \lim \sup_{u\downarrow 0\text{ in }\mathbb{A}}F^{\ast
}(v+u)\leq \lim \sup_{u\downarrow 0}F^{\ast }(v+u)\leq F^{\ast }(v)+F^{\ast
}(0+)=F^{\ast }(v),
\]%
and so%
\[
F^{\ast }(v+)=F^{\ast }(v),
\]%
i.e. $F^{\ast }$ is right-continuous everywhere on $\mathbb{R}_{+}$. Now for
$u\in \mathbb{A}$ with $0<u<w$%
\[
F^{\ast }(w-u)<F^{\ast }((w-u)+u)=F^{\ast }(w).
\]%
So, for arbitrary $0<v<w,$ and $u\in \mathbb{A}$ with $u>0$ such that $%
v<w-u<w,$%
\[
F^{\ast }(v)=F^{\ast }(v+)=\lim \inf \{F^{\ast }(w-u):v<w-u<w,u\in \mathbb{A}%
\}<F^{\ast }(w),
\]%
as $\mathbb{A}$ is dense. So%
\[
F^{\ast }(v)<F^{\ast }(w),
\]%
i.e. $F^{\ast }$ is increasing. So it is continuous at some $v>0.$ Since $%
K(0)=0,$ $K$ is continuous at $0$ on $\mathbb{A}$. But%
\[
F^{\ast }(v-u)-F^{\ast }(v)=g(v)K(-u)\text{\qquad }(u\in \mathbb{A})(v\in
\mathbb{R}),
\]%
so $F^{\ast }$ satisfies, for any $v\in \mathbb{R}_{+},$%
\[
\lim_{u\downarrow 0\text{ \& }u\in \mathbb{A}}F^{\ast }(v-u)=F^{\ast }(v).
\]%
But $F^{\ast }$ is increasing, and $\mathbb{A}$ is dense; so $F^{\ast }$ is
continuous. Then for any fixed $u>0$ in $\mathbb{A},$%
\[
g(v)=[F^{\ast }(v+u)-F^{\ast }(v)]/K(u),
\]%
and so $g$ is continuous at any $v\in \mathbb{R}_{+},$ since $F^{\ast }$ is
continuous at $v$ and at $v+u.$ $\square $

\section{From the Goldie to the Beurling Equation}

In $(GFE)$, take $K$ and $F^{\ast }$ (which will reduce to the same -- see
Th. 5 below) as $K.$ We generalize the $e^{-\rho \cdot }$ to $g$, which will
serve as an \textit{auxiliary function }(which will reduce to $e^{-\rho
\cdot }$ in the case of interest). We now have the Goldie equation in the
form%
\[
K(v+u)-K(v)=g(v)K(u).
\]%
For reasons that will emerge (see inter alia \S 5), an important
generalization arises if on the left the additive action of $v$ on $u$ is
made dependent on $g$:
\begin{equation}
K(v+ug(v))-K(v)=g(v)K(u),  \label{double-star}
\end{equation}%
so that while $g$ appears twice, $K$ still appears here three times. This
form is closely related to a situation with all function symbols identical, $%
\varphi $ say (which we will take non-negative):
\begin{equation}
\varphi (v+u\varphi (v))=\varphi (u)\varphi (v),\text{ \qquad }(\forall
u,v\in \mathbb{R}_{+}).  \tag{$BFE$}
\end{equation}%
Indeed, from here, writing $g$ for $\varphi $ and with $K(t)\equiv g(t)-1$
(i.e. as in (*) with $\kappa =1$), we recover (\ref{double-star}).

This $(BFE)$ is our \textit{Beurling functional equation, }a special case of
the Go\l \k{a}b-Schinzel equation (see \S 1) in view of the non-negativity
and of the domain being $\mathbb{R}_{+}$ rather than $\mathbb{R}$ (both
considerations arising from the context of Beurling regular variation).
Solutions of the `conditional' Go\l \k{a}b-Schinzel equation (i.e. with
domain restricted to $\mathbb{R}_{+}$, but without the non-negativity
restriction) were considered and characterized in [BrzM] and shown to be
extendible uniquely to solutions with domain $\mathbb{R}.$ Note that for any
extension to $\mathbb{R}_{+}\cup \{0\},$ if $\varphi (0)=0,$ then $(BFE)$
implies $\varphi =0;$ we will therefore usually set $\varphi (0)=1,$ the
alternative dictated by the equation $\varphi (0)=\varphi (0)^{2}.$
Solutions $\varphi >0$ are relevant to the Beurling theory of regular
variation -- see [Ost2] for an analysis; their study is much simplified by
the following easy result, inspired by a close reading of [Brz1, Prop. 2].

\bigskip

\noindent \textbf{Theorem 5.} \textit{If} $\varphi :\mathbb{R}%
_{+}\rightarrow \mathbb{R}_{+}$ \textit{satisfies }$(BFE),$\textit{\ then }$%
\varphi (x)\geq 1$\textit{\ for all }$x>0.$

\bigskip

\noindent \textit{Proof.} Suppose that $\varphi (u)<1$ for some $u>0;$ then $%
v:=u/(1-\varphi (u))>0$ and so, since $v=u+v\varphi (u),$%
\[
0<\varphi (v)=\varphi \left( u+v\varphi (u)\right) =\varphi (u)\varphi (v).
\]%
So cancelling $\varphi (v),$ one has $\varphi (u)=1,$ a contradiction. $%
\square $

\bigskip

The theorem above motivates the introduction of an important tool in the
study of positive solutions $\varphi $: the \textit{Beck sequence} $%
t_{m}=t_{m}(u)$, defined for any $u>0$ recursively by
\[
t_{m+1}=t_{m}+u\varphi (t_{m})\text{ with }t_{0}=0,
\]%
so that%
\[
\varphi (t_{m+1})=\varphi (u)\varphi (t_{m}).
\]%
By Theorem 5, the sequence $t_{m}$ is divergent, since either $\varphi (u)=1$
and $t_{m}=mu,$ or else%
\begin{equation}
t_{m}=u\frac{\varphi (u)^{m}-1}{\varphi (u)-1}=(\varphi (u)^{m}-1)\left/
\frac{\varphi (u)-1}{u}\right. ,  \label{formula}
\end{equation}%
e.g. by Lemma 4 of [Ost2] (cf. a lemma of Bloom: BGT Lemma 2.11.2). In
either case, for $u,t>0$ a unique integer $m=m_{t}(u)$ exists satisfying%
\[
t_{m}\leq t<t_{m+1}.
\]%
This tool will enable us to prove in Theorem 7 below that a positive
solution of $(BFE)$ takes the form $\varphi (t)=1+ct$ for some $c\geq 0.$
Theorem 6 and its Corollary below lay the foundations.

\bigskip

\noindent \textbf{Theorem 6.} \textit{If a function }$\varphi \geq 0$\textit{%
\ satisfies the equation }($BFE$)\textit{\ on }$\mathbb{R}_{+}$ \textit{with
}$\varphi (t)>1$\textit{\ for }$t\in I=(0,\delta ),$\textit{\ for some }$%
\delta >0,$ \textit{then }$\varphi $ \textit{is continuous and (strictly)
increasing, and }$\varphi >1$\textit{.}

\bigskip

\noindent \textit{Proof}\textbf{.} Take $K(t)=\varphi (t)-1;$ then $K>0$ on $%
I.$ Writing $x=u$ and $y=v\varphi (u),$%
\[
\varphi (x+y)-\varphi (x)=K(y/\varphi (x))\varphi (x).
\]%
Fix $x\in I$; then $\varphi (x)>1,$ and so $y/\varphi (x)\in I$ for $y\in I,$
so that $K(y/\varphi (x))>0$. As in Theorem 2, $\varphi (x+y)>\varphi (y)$
for $x,y\in I,$ and $\varphi $ is increasing on a subinterval of $I$. So $%
\varphi $ is continuous at some point $u\in I,$ $\varphi (u)>0$ and
\[
\varphi (u)=\lim_{v\downarrow 0}\varphi (u+v\varphi (u))=\varphi
(u)\lim_{v\downarrow 0}\varphi (v):\qquad \varphi (0+)=\lim_{v\downarrow
0}\varphi (v)=1.
\]%
So for $x>0$ with $\varphi (x)>0,$%
\[
\lim_{v\downarrow 0}\varphi (x+v\varphi (x))=\varphi (x)\lim_{v\downarrow
0}\varphi (v)=\varphi (x),
\]%
and so $\varphi $ is right-continuous at $x.$

Let $J\supseteq I$ be a maximal interval $(0,\eta )$ on which $\varphi $ is
increasing, and suppose that $\eta $ is finite. Then $\varphi (\eta )=0:$
otherwise, $\varphi (\eta )>0$ and as above $\varphi $ is right-continuous
at $\eta ;$ since $\varphi >1$ in $I,$ $\varphi (t)>1$ near $\eta ;$ then $%
\varphi $ is increasing to the right of $\eta ,$ a contradiction. Now choose
$t<\eta <t+\delta ;$ then $v=(\eta -t)/\varphi (t)<x-t,$ and%
\[
\varphi (\eta )=\varphi (t+v\varphi (t))=\varphi (t)\varphi (v)>0,
\]%
as $t\in J,$ contradicting $\varphi (\eta )=0.$ This shows that $J=\mathbb{R}%
_{+}$, and so $\varphi $ is right-continuous and increasing.

Finally, we check that $\varphi $ is left-continuous at any $x>0.$ Let $%
z_{n}\downarrow 0$ with $x-z_{n}>0;$ then, as above, $\varphi
(z_{n})\rightarrow 1.$ As $\varphi $ is increasing, $u_{n}:=\varphi
(x-z_{n})-1\leq \varphi (x)-1$ is positive and bounded, so%
\begin{eqnarray*}
\varphi (x-z_{n}) &=&\varphi ((x-z_{n})+z_{n}\varphi (x-z_{n}))/\varphi
(z_{n}) \\
&=&\varphi (z+z_{n}u_{n})/\varphi (z_{n})\rightarrow \varphi (x),
\end{eqnarray*}%
by right-continuity at $x$. $\square $

\bigskip

\noindent \textbf{Corollary.} \textit{If }$\varphi >0,$\textit{\ then }$%
\varphi $\textit{\ is continuous, and either }$\varphi >1,$\textit{\ or the
value 1 is repeated densely and so }$\varphi \equiv 1.$

\bigskip

\noindent \textit{Proof. }By Th. 5 $\varphi \geq 1,$ so $\varphi $ is
(weakly) increasing and so continuous. Suppose that $\varphi >1$ is false;
then, by Theorem 6, there is no interval $(0,\delta )$ with $\delta >0$ on
which $\varphi >1,$ and so there are arbitrarily small $u>0$ with $\varphi
(u)=1.$ Fix $t>0.$ For any $u$ with $\varphi (u)=1,$ choose $n=n_{t}(u)$
with $t_{n}:=nu\leq t\leq (n+1)u,$ as above. Then $\varphi (t_{n})=1$ and $%
0\leq t-t_{n}<u.$ So the value $1$ is taken on a dense set of points, and so
by continuity $\varphi (t)\equiv 1.$ $\square $

\bigskip

We now adapt Goldie's argument above to give an easy proof of the following.
Theorem 7 below can be derived from [Brz1, Cor 3] or [Brz2, Th1]. There
algebraic considerations are key; an analytical proof was provided in
[Ost2], but by a different and more complicated route. We include the proof
below for completeness, as it is analogous to the Goldie Theorem above and
so thematic here. We use a little less than Theorem 6 provides.

\bigskip

\noindent \textbf{Theorem 7.} \textit{If }$\varphi (t)>1$ \textit{holds for
all }$t$\textit{\ in some interval }$(0,\delta )$ \textit{with }$\delta >0,$
\textit{and satisfies }$(BFE)$\textit{\ on }$\mathbb{R}_{+}$\textit{, then }$%
\varphi $\textit{\ is differentiable, and takes the form}
\[
\varphi (t)=1+ct.
\]%
\noindent \textit{Proof.} Fix $x_{0}>0$ with $\varphi (x_{0})\neq 1.$ Put%
\[
K(t):=\varphi (t)-1.
\]%
By Th. 6 $K$ is continuous, so $K(t)\neq 0$ for $t$ sufficiently close to $%
x_{0};$ we may assume also that $\varphi (u)\neq 1$ for all small enough $%
u>0,$ and so $K(u)\neq 0$ for sufficently small $u.$

Let $x$ be arbitrary; in the analysis below $x$ and $x_{0}$ play similar
roles, so it will be convenient to also write $x_{1}$ for $x.$

For $j=0,1$ and any $u>0,$ referring to the Beck sequence $t_{m}=t_{m}(u)$
as above, select $i_{j}=i_{j}(u):=m_{x_{j}}(u)$ so that%
\[
t_{i_{j}}\leq x_{j}<t_{i_{j}+1}\qquad (j=0,1);
\]%
then%
\[
\varphi (t_{m+1})-\varphi (t_{m})=\varphi (u)\varphi (t_{m})-\varphi
(t_{m})=K(u)\varphi (t_{m}).
\]%
Summing,%
\[
\varphi (v_{m})-\varphi (v_{0})=K(u)\sum_{n=0}^{m-1}\varphi (t_{n}).
\]%
As noted, for all small enough $u$, $K(v(i_{0}))$ non-zero (this use
compactness of $[0,x_{0}])$. Cancelling $K(u)$ below (as also $K(u)$ is
non-zero), and introducing $u$ in its place (to get the telescoping sums),
\[
\frac{K(t_{i_{1}})}{K(t_{i_{0}})}=\frac{K(u)\sum_{n=0}^{i_{1}-1}\varphi
(t_{n})}{K(u)\sum_{n=0}^{m-1}\varphi (t_{n})}=\frac{\sum_{n=0}^{i_{1}-1}u%
\varphi (t_{n})}{\sum_{n=0}^{m-1}u\varphi (t_{n})}=\frac{%
\sum_{n=0}^{i_{1}-1}t_{n+1}-t_{n}}{\sum_{n=0}^{i_{0}-1}t_{n+1}-t_{n}}=\frac{%
t_{i_{1}}}{t_{i_{0}}},
\]%
as%
\[
t_{n+1}-t_{n}=u\varphi (t_{n}).
\]%
Passing to the limit as $u\rightarrow 0,$ by continuity%
\[
K(x_{1})/K(x_{0})=x_{1}/x_{0}.
\]%
Setting $c_{0}:=K(x_{0})/x_{0},$
\[
\varphi (x)-1=K(x)=c_{0}x:\text{ }\varphi (x)=1+c_{0}x.\text{ }\square
\]

\bigskip

\noindent \textit{Remark. }The argument above could have been presented
explicitly in terms of integrals. Relevant here is the differentiability
referred to above.

\section{On the Beurling functional equation}

This section is suggested by the recent work of the second author [Ost2]. We
include it here for two reasons. First, it is thematically close to other
results in this paper. Secondly, it provides a simpler proof of results
concerning Beurling's equation than through specialization of results by Brzd%
\k{e}k and by Brzd\k{e}k and Mure\'{n}ko in [Brz1] and [BrzM]. Theorem B at
the end of this section is taken from these papers; we include it here for
completeness, and as our proof is more direct and shorter.

In Th. 8 below, with context as in Ths. 5 and 6, we use $\lambda $ rather
than $\varphi $ for ease of comparison with [Ost2]. For $\lambda :[0,\infty
)\rightarrow \mathbb{R},$ denote its level set above unity by:%
\[
L_{+}(\lambda ):=\{t\in \mathbb{R}_{+}:\lambda (t)>1\}.
\]

\bigskip

\noindent \textbf{Theorem 8. }\textit{If the continuous solution }$\lambda $
\textit{of }$(BFE)$\textit{\ with }$\lambda (0)=1$ \textit{has a nonempty
level set }$L_{+}(\lambda )$ \textit{containing an interval }$(0,\delta )$
\textit{for some }$\delta >0$\textit{, then }$\lambda $ \textit{is
differentiable and for some }$\rho >0$\textit{\ }%
\[
\lambda (t)\equiv 1+\rho t.
\]%
\textit{\ }

\noindent \textit{Proof.} We recall (from the proof of Theorem 5 above) the
Beck sequence $t_{n}(u),$ and that $\lambda (t_{n}(u))=\lambda (u)^{n};$
from here, for $u\in L_{+}$, by summing,%
\[
t_{n}(u)=u\frac{\lambda (u)^{n}-1}{\lambda (u)-1}=(\lambda (u)^{n}-1)\left/
\frac{\lambda (u)-1}{u}\right. ,
\]%
(and with $t_{n}(u)=nu$ if $\lambda (u)=1),$ and from the recurrence,%
\[
\Delta _{m}(u):=t_{m+1}(u)-t_{m}(u)=u\lambda (u)^{m}.
\]%
For $T\in L_{+}$ and $u>0,$ write $m=m(u)=m_{T}(u)$ for the jump index for
which%
\[
t_{m}(u)\leq T<t_{m+1}(u).
\]%
By (\ref{formula}) and continuity at $0$ of $\lambda $,
\begin{equation}
\Delta _{m(u)}(u)=u\lambda (u)^{m(u)}\leq T(\lambda (u)-1)+u\rightarrow 0%
\text{ as }u\rightarrow 0,  \label{del}
\end{equation}%
for $u\in L_{+}$ \textit{uniformly} in $T>0$ on compacts. Likewise for $%
u\notin L_{+}$, as then $\Delta _{m(u)}(u)=u.$

Consider any null sequence $u_{n}\rightarrow 0$ with $u_{n}>0.$ We will show
that $\{(\lambda (u_{n})-1)\left/ u_{n}\right. \}$ is convergent, by showing
that down every subsequence $\{(\lambda (u_{n})-1)\left/ u_{n}\right.
\}_{n\in \mathbb{M}}$ there is a convergent sub-subsequence with limit
independent of $\mathbb{M}$.

W.l.o.g. we take $0<u_{n}\in L_{+}$ for all $n$ (so $u_{n}<\delta ).$ Now
consider an arbitrary $T\in L_{+}.$ Passing to a subsequence (dependent on $%
T)$ of $\{(\lambda (u_{n})-1)\left/ u_{n}\right. \}_{n\in \mathbb{M}}$ if
necessary, we may suppose, for $k(n):=m_{T}(u_{n}),$ that $\Delta
_{k(n)}(u_{n})\rightarrow 0;$ then along $\mathbb{M}$
\[
|T-t_{m(u_{n})}(u_{n})|\leq \Delta _{m(u_{n})}(u_{n}),\text{ and so }%
t_{k(n)}(u_{n})=t_{m(u_{n})}(u_{n})\rightarrow T.
\]%
Again by (\ref{formula}) and continuity at $T$ of $\lambda ,$ putting $\rho
:=(\lambda (T)-1)/T>0,$%
\[
\frac{\lambda (u_{n})-1}{u_{n}}=\frac{\lambda (u_{n})^{m(u_{n})}-1}{%
t_{m(u_{n})}(u_{n})}=\frac{\lambda (t_{m(n)}(u_{n}))-1}{t_{m(u_{n})}(u_{n})}%
\rightarrow \frac{\lambda (T)-1}{T}=\rho ,
\]%
along $\mathbb{M}$ to a limit $\rho $ dependent only on $T.$ So $\{(\lambda
(u_{n})-1)\left/ u_{n}\right. \}$ is itself convergent to $\rho .$ But this
holds for any null sequence $\{u_{n}\}$ in $\mathbb{R}_{+},$ so the function
$\lambda $ is differentiable at $0,$ and so is right-differentiable
everywhere in $L_{+}$ (see [Ost2, Lemma 3]). It is also left-differentiable
at any $x>0,$ as follows. For $y$ with $0<y<x,$ put $t:=(x-y)/\lambda (y)>0.$
Then $x=y+t\lambda (y),$ so%
\[
\frac{\lambda (x)-\lambda (y)}{x-y}=\frac{\lambda (y+t\lambda (y))-\lambda
(y)}{x-y}=\frac{[\lambda (t)-1]\lambda (y)}{x-y}=\frac{\lambda (t)-1}{t}.
\]%
But $t\downarrow 0$ as $y\uparrow x$ (by continuity of $\lambda $ at $x),$
and $(\lambda (t)-1)/t\rightarrow \lambda ^{\prime }(0).$ So $\lambda $ is
left-differentiable at $x$ and so differentiable; from here $\lambda
^{\prime }(x)=\lambda ^{\prime }(0).$

Integration then yields $\lambda (x);$ also, since $T$ above was arbitrary,
for any $T\in L_{+}$
\[
\rho =\lim_{_{n\in \mathbb{M}}}\{(\lambda (u_{n})-1)\left/ u_{n}\right.
\}=\lambda ^{\prime }(0)=(\lambda (T)-1)/T:\qquad \lambda (x)=1+\rho x\text{
}(x\in \mathbb{R}_{+}).\text{ }\square
\]

In Theorem BM below we use $f$ rather than $\varphi $ for ease of comparison
with [BrzM].

\bigskip

\noindent \textbf{Theorem BM }([BrzM, Lemma 7])\textbf{. }\textit{For }$f>0$%
\textit{\ a solution of }$(BFE)$\textit{, if }$f\neq 1$\textit{\ at all
points, then }$f(x)=1+cx$ \textit{for some }$c>0.$

\bigskip

\noindent \textit{Proof. }By symmetry, for any $x,y$%
\[
f(x+yf(x))=f(x)f(y)=f(y+xf(y)).
\]%
Fix $x$ and $y$ and put $u=x+yf(x)$ and $v=y+xf(y).$ If these are unequal,
w.l.o.g. suppose that $v>u.$ Then $(v-u)/f(u)>0,$ so%
\[
0<f(u)=f(v)=f(u+f(u)(u-v)/f(u))=f(u)f((u-v)/f(u)).
\]%
So $f((u-v)/f(u))=1,$ a contradiction. So $u=v$: that is, for all $x,y>0$%
\[
x+yf(x)=y+xf(y);
\]%
equivalently, for all $x,y>0$%
\[
x/(1-f(x))=y/(1-f(y))=c,
\]%
say. Then $f(x)=1+cx$ for all $x>0.$ So $c>0.$ $\square $

\bigskip

Below we suppose that $f(a)=1,$ for some fixed $a>0.$ Note that $t_{n}:=na$
is a \textit{Beck} sequence with step size $a;$ so $f(na)=1,$ since $%
f(t_{n})=f(t_{1})^{n}$.

\bigskip

For $f>0$ a solution of $(BFE),$ we denote here the range of $f$ by $%
R_{f}:=\{w:(\exists x>0)w=f(x)\}.$ If $f\equiv 1,$ then $R_{f}=\{1\}.$

\bigskip

\noindent \textbf{Lemma B }([Brz1, Cor.1], cf. [BrzM, Lemmas 1,2])\textbf{.}
\textit{If the value }$1$\textit{\ is achieved by a solution }$f>0$\textit{\
of }$(BFE),$ \textit{then}

\noindent (i) \textit{the range set }$R_{f}$ \textit{is a multiplicative
subgroup;}

\noindent (ii)\textbf{\ }$f(x+a)=f(x)$ for all $x>0;$

\noindent (iii) $f(wa)=1$ for $w\in R_{f}.$

\bigskip

\noindent \textit{Proof.} For (i), $(BFE)$ itself implies that $R_{f}$ is a
semigroup. We only need to find the inverse of $w:=f(x)$ with $x>0.$ Choose $%
n\in \mathbb{N}$ with $na>x.$ Put $y=(na-x)/f(x);$ then%
\[
f(x)f(y)=f(x+yf(x))=f(na)=1.
\]%
For (ii), note that $f(x)f(a)=f(a+xf(a))=f(x+a).$ For (iii), this time write
$w=1/f(x);$ then%
\[
f(x)=f(x+a)=f(x+f(x)a/f(x))=f(x)f(aw),
\]%
and cancelling $f(x)>0$ gives $f(aw)=1.$ $\square $

\bigskip

\noindent \textbf{Theorem B }([Brz1, Th. 3])\textbf{.} \textit{If }$1\in
R_{f}$\textit{, then} $f\equiv 1.$

\bigskip

\noindent \textit{Proof.} Suppose otherwise; then, by Theorem 5, $f(u)>1,$
for some $u>0.$ Choose $n\in \mathbb{N}$ with $na>u/(f(u)-1)>0,$ and put%
\[
v:=na+u/(1-f(u))>0:\text{ }v+naf(u)=na+u+vf(u).
\]%
So, since $f(u)\in R_{f},$ applying Lemma B (first (ii), then (i))
\begin{eqnarray*}
0 &<&f(v)=f(v+f(u)na)=f(u+vf(u)+na) \\
&=&f(u+vf(u))=f(u)f(v),
\end{eqnarray*}%
yielding the contradiction $f(u)=1.$ Hence $f(x)=1$ for all $x.$ $\square $

\section{Extensions of the Goldie and Beurling equations}

Below we consider two generalizations of the Beurling equation inspired by
Goldie's equation, relevant to Beurling regular variation, for which see
[BinO3]. The first uses three functions:
\begin{equation}
K(v+uk(v))-K(v)=g(u)k(v)\qquad (u,v\in \mathbb{R}_{+})  \tag{$GBE$}
\end{equation}%
Here the choice $k=K=\varphi $ with $g=\varphi (u)-1$ recovers the Beurling
equation. One can also form a Pexider-like generalization (for which see
[Kuc,13.3], or [AczD, 4.3]) for the right-hand side above, by replacing the
occurrence of $k$ there with an additional function $h:$
\begin{equation}
K(v+uk(v))-K(v)=g(u)h(v).\qquad (u,v\in \mathbb{R}_{+})  \tag{$GBE$-$P$}
\end{equation}%
Here $h=K$ and $k=1$ yields Goldie's equation; $h=k=K,g=1-k$ yields the
Beurling equation.

\bigskip

\noindent \textbf{Theorem 9.} \textit{Consider the functional equation }$%
(GBE $-$P)$\textit{\ with }$g(u)\neq 0$\textit{\ for }$u>0$ \textit{near }$%
0, $ $h,k$\textit{\ continuous on }$\mathbb{R}_{+}\cup \{0\}$\textit{\ and
positive,\ and with }$k(0)>0.$\newline
\textit{With}%
\[
H(x):=\int_{0}^{x}h(t)\frac{dt}{k(t)},\text{ for }x\geq 0,
\]%
\textit{any solution }$K$ \textit{is differentiable and takes the form }$%
K(x)=cH(x)$\textit{, for }$c$\textit{\ a constant; furthermore, }$g$ \textit{%
is continuous with }$g(0+)=0$ \textit{and }$h(0)=0.$

\bigskip

\noindent \textit{Proof}\textbf{.} Since%
\[
K(v+w)-K(v)=g(w/k(v))h(v),
\]%
we deduce for $u,v\in \mathbb{R}_{+}$ that%
\[
K(v+)-K(v)=g(0+)h(v),\qquad \text{ }K(u-)-K(u)=g(0+)h(u),
\]%
and that for any $v>0$ and all small enough $w>0$%
\[
K(v+w)>K(v).
\]%
So $K$ is locally increasing on $\mathbb{R}_{+}$, and so is increasing and
so continuous on a dense set $D\subseteq \mathbb{R}_{+}$. For $v\in D,$%
\[
g(0+)h(v)=K(v+)-K(v)=0,
\]%
and since $h>0,$ $g(0+)=0.$ So $K$ is continuous on $\mathbb{R}_{+}$.

Now consider for $u>0$ the Beck sequence
\[
t_{n+1}(u)=t_{n}(u)+uk(t_{n}(u)),\text{ }t_{0}=0,
\]%
which is increasing as $k>0.$ For any $t,u>0$ we claim there is $m=m_{t}(u)$
with%
\[
m_{t}(u)\leq t<m_{t}(u)+1.
\]%
For otherwise, with $t,u$ fixed as above, the sequence $t_{n}(u)$ is bounded
by $t$ and, putting $\tau :=\sup t_{n}(u)\leq t$,
\[
k(t_{n}(u)=\frac{1}{u}[t_{n+1}(u)-t_{n}(u)]\rightarrow 0,
\]%
contradicting lower boundedness of $k$ near $\tau .$ Next observe that,
since $k$ is bounded on $[0,t],$ by $M_{t}$ say,%
\[
t_{m+1}(u)-t_{m}(u)=uk(t_{k}(u)\leq uM_{t}\rightarrow 0.
\]

Now fix $x_{0,}x_{1}>0.$ Select $i_{0}=i_{0}(u)$ and $i_{1}=i_{1}(u)$ so that%
\[
t_{i_{j}}\leq x_{j}<t_{i_{j}+1}.
\]%
Then%
\[
K(t_{m+1})-K(t_{m})=g(u)h(t_{m}).
\]%
Summing, and setting $p(t):=h(t)/k(t),$%
\[
K(t_{m})-K(t_{0})=g(u)\sum_{n=0}^{m-1}h(t_{n})=\frac{g(u)}{u}%
\sum_{n=0}^{m-1}uk(t_{n})p(t_{n}).
\]%
For all small enough $u$ we have $g(u)$ non-zero, so
\[
\frac{K(t_{i_{1}})}{K(t_{i_{0}})}=\frac{g(u)\sum_{n=0}^{i_{1}-1}h(t_{n})}{%
g(u)\sum_{n=0}^{m-1}h(t_{n})}=\frac{\sum_{n=0}^{i_{1}-1}uk(t_{n})p(t_{n})}{%
\sum_{n=0}^{m-1}uk(t_{n})p(t_{n})}\rightarrow \frac{\int_{0}^{x_{i}}p(t)dt}{%
\int_{0}^{x_{0}}p(t)dt}.
\]%
Passing to the limit as $u\rightarrow 0,$ by continuity of $K,$%
\[
K(x_{1})/K(x_{0})=H(x_{1})/H(x_{0}).
\]%
Setting $c_{0}:=K(x_{0})/H(x_{0}),$ we have%
\[
K(x)=c_{0}H(x),
\]%
for $x\geq 0,$ which is differentiable.

Substitution in $(GBE$-$P)$ gives%
\[
g(u)=\frac{1}{h(v)}\int_{v}^{v+uk(v)}h(t)\frac{dt}{k(t)},
\]%
which for any $v$ is continuous in $u;$ further,%
\[
g(u)=\left( 1/h(v)\right) \cdot \left( h(v+u\theta k(v))/k(v+u\theta
k(v))\right) ,
\]%
for some $\theta =\theta _{v}(u)$ with $0<\theta _{v}(u)<1.$ Assuming $%
h(0)\neq 0,$ taking limits as $v\rightarrow 0+$ we have for some $0\leq
\theta _{0}(u)\leq 1$ that%
\[
g(u)=\left( 1/h(0)\right) \cdot \left( h(u\theta _{0}k(0))/k(u\theta
_{0}k(0))\right) .
\]%
Now take limits as $u\rightarrow 0:$ as $g(0+)=0,$%
\[
0=g(0)=1/k(0)>0,
\]%
a contradiction. So $h(0)=0.$ $\square $

\bigskip

Taking $h=k$ so that $p=1,$ which is already continuous, we obtain a
corollary which needs only local boundedness above and away from 0, rather
than continuity in $k$ (to justify the use of the Beck sequence).

\bigskip

\noindent \textbf{Theorem 10.} \textit{Consider the functional equation }$%
(GBE)$\textit{\ above with }$k>0$\textit{\ locally bounded above and away
from }$0$\textit{\ on }$\mathbb{R}_{+}$\textit{, and }$g(u)\neq 0$\textit{\
for }$u\neq 0$ \textit{near }$0$\textit{.}

\textit{Any solution is linear: }$K(x)=cx$\textit{; furthermore, }$g(u)=cu$.
\textit{In particular, for }$K=k$ \textit{and }$g=k-1,\mathit{\ }$\textit{%
the solution of the Beurling equation is }$k(u)=1+cu.$

\bigskip

\noindent \textit{Proof}\textbf{.} This is a simpler version of the
preceding proof. With the notation there and with $h=k$ (so that $p=1$ and $%
t_{n+1}-t_{n}=uk(t_{n})),$ the proof reduces to noting that%
\[
K(t_{m})-K(t_{0})=g(u)\sum_{n=0}^{m-1}k(t_{n})=\frac{g(u)}{u}%
\sum_{n=0}^{m-1}uk(t_{n})=\frac{g(u)}{u}\sum_{n=0}^{i_{1}-1}(t_{n+1}-t_{n})=%
\frac{g(u)}{u}t_{i_{1}},
\]%
and that for all small enough $u$ we have $g(u)$ non-zero. Then
\[
\frac{K(t_{i_{1}})}{K(t_{i_{0}})}=\frac{g(u)\sum_{n=0}^{i_{1}-1}k(t_{n})}{%
g(u)\sum_{n=0}^{m-1}k(t_{n})}=\frac{\sum_{n=0}^{i_{1}-1}uk(t_{n})}{%
\sum_{n=0}^{m-1}uk(t_{n})}=\frac{\sum_{n=0}^{i_{1}-1}t_{n+1}-t_{n}}{%
\sum_{n=0}^{i_{0}-1}t_{n+1}-t_{n}}=\frac{t_{i_{1}}}{t_{i_{0}}},
\]%
Passing to the limit as $u\rightarrow 0,$ by continuity of $K,$%
\[
K(x_{1})/K(x_{0})=x_{1}/x_{0}:\qquad K(x)=c_{0}x,
\]%
with $c_{0}:=K(x_{0})/x_{0}$. $\square $

\section{Complements}

\noindent 1. \textit{Regular variation: related results. }Our results here
concern the `Goldie argument', the crux of the remaining `hard proof' in
regular variation (see the proof of Theorem 3 above and [BinO4, \S 8]). We
have focussed particularly here on the key relevant results in BGT, namely
Theorem 1.4.3 and Th. 3.2.5, the former simplified in [BinO4, Th. 6], the
latter here in Theorem 3. There are a number of related and similar results
in BGT Ch. 3. and these may be treated similarly.

\noindent 2. \textit{Go\l \k{a}b-Schinzel and related functional equations.}
For a recent text-book account of the equation see [AczD, Ch. 19] or the
more recent surveys [Brz3] or [Jab], which include generalizations and a
discussion of applications in algebra, meteorology and fluid mechanics --
see for instance [KahM].

\noindent 3. \textit{Symmetrized Goldie functional equations.} The equation%
\[
K(x+y)=g(x)K(y)+g(y)K(x)
\]%
is studied in [AczD, Ch. 13] in connection with trigonometric identities;
interpreting the right-hand side as a discrete convolution with $g$ positive
and $g(x)+g(y)=1$ here also connects to the context of Markov chains -- for
which see [AczD], Ch. 12]. The `Goldie case' $g(x)=e^{-\rho x}/2$ yields $%
K(x)=e^{-\rho x}.$

\noindent 4. \textit{Differentiability from continuity.} Various equations
are known to confer additional regularity on solutions. Here the classical
example is the Euler-Lagrange equation; another is provided by d'Alembert's
wave equation, for which see [AczD, Ch. 14].

\noindent 5. \textit{Flows.} The integrator $du/k(u)$ in Theorem 9 above is
connected to the \textit{flow} aspects of Beurling regular variation, for
which see [BinO3], [Ost1]. The background involves topological dynamics,
inspired by Beck [Bec, 5.41]. It also involves regular variation in more
general contexts than $\mathbb{R}$, such as normed groups; see [BinO1] and
the references cited there for detail.

\bigskip

\begin{center}
\textbf{References}
\end{center}

\noindent \lbrack AczD] J. Acz\'{e}l, J. Dhombres, \textsl{Functional
equations in several variables. With applications to mathematics,
information theory and to the natural and social sciences.} Encyclopedia of
Math. and its App., 31, CUP, 1989\newline
\noindent \lbrack Bec] A. Beck, \textsl{Continuous flows on the plane},
Grundl. math. Wiss. \textbf{201}, Springer, 1974. \newline
\noindent \lbrack Bin] N. H. Bingham, Tauberian theorems and the central
limit theorem, \textsl{Ann. Prob.} \textbf{9} (1981), 421-431.\newline
\noindent \lbrack BinG] N. H. Bingham, C. M. Goldie, Extensions of regular
variation: I. Uniformity and quantifiers, \textsl{Proc. London Math. Soc.}
(3) \textbf{44} (1982), 473-496.\newline
\noindent \lbrack BinGT] N. H. Bingham, C. M. Goldie and J. L. Teugels,
\textsl{Regular variation}, 2nd ed., Cambridge University Press, 1989 (1st
ed. 1987). \newline
\noindent \lbrack BinO1] N. H. Bingham and A. J. Ostaszewski, Normed groups:
Dichotomy and duality. \textsl{Dissertationes Math.} \textbf{472} (2010),
138p. \newline
\noindent \lbrack BinO2] N. H. Bingham and A. J. Ostaszewski, Dichotomy and
infinite combinatorics: the theorems of Steinhaus and Ostrowski. \textsl{%
Math. Proc. Cambridge Phil. Soc.} \textbf{150} (2011), 1-22. \newline
\noindent \lbrack BinO3] N. H. Bingham and A. J. Ostaszewski, Beurling slow
and regular variation, \textsl{Trans. London Math. Soc., }to appear (see
also arXiv1301.5894 and arXiv 1307.5305).\newline
\noindent \lbrack BinO4] N. H. Bingham and A. J. Ostaszewski: Additivity,
subadditivity and linearity: automatic continuity and quantifier weakening,
arXiv.1405.3948.\newline
\noindent \lbrack BojK] R. Bojani\'{c} and J. Karamata, \textsl{On a class
of functions of regular asymptotic behavior, }Math. Research Center Tech.
Report 436, Madison, Wis. 1963; reprinted in \textsl{Selected papers of
Jovan Karamata} (ed. V. Mari\'{c}, Zevod Ud\v{z}benika, Beograd, 2009),
545-569.\newline
\noindent \lbrack Brz1] J. Brzd\k{e}k, Subgroups of the group $\mathbb{Z}%
_{n} $ and a generalization of the Go\l \k{a}b-Schinzel functional equation,
\textsl{Aequat. Math.} \textbf{43} (1992), 59-71.\newline
\noindent \lbrack Brz2] J. Brzd\k{e}k, Bounded solutions of the Go\l \k{a}%
b-Schinzel equation, \textsl{Aequat. Math.} \textbf{59} (2000), 248-254.%
\newline
\noindent \lbrack Brz3] J. Brzd\k{e}k, The Go\l \k{a}b-Schinzel equation and
its generalizations, \textsl{Aequat. Math.} \textbf{70} (2005), 14-24.%
\newline
\noindent \lbrack BrzM] J. Brzd\k{e}k and A. Mure\'{n}ko, On a conditional Go%
\l \k{a}b-Schinzel equation, \textsl{Arch. Math.} \textbf{84} (2005),
503-511.\newline
\noindent \lbrack Dal1] H. G. Dales, \textsl{Automatic continuity: a survey}%
, Bull. London Math. Soc. \textbf{10} (1978),129-183.\newline
\noindent \lbrack Dal2] H. G. Dales, \textsl{Banach algebras and automatic
continuity}, London Math. Soc. Monographs 24, Oxford University Press, 2000.%
\newline
\noindent \lbrack GolS] St. Go\l \k{a}b and A. Schinzel, Sur l'\'{e}quation
fonctionelle, $f[x+yf(x)]=f(x)f(y),$ \textsl{Publ. Math. Debrecen}, \textbf{6%
} (1959), 113-125.\newline
\noindent \lbrack Ger] R. Ger, Almost additive functions onsemigroups and a
functional equation, \textsl{Publ. Math. Debrecen} \textbf{26} (1979),
219-228.\newline
\noindent \lbrack Hei] C. H. Heiberg, A proof of a conjecture by Karamata.
\textsl{Publ. Inst. Math. (Beograd)} (N.S.) \textbf{12} (26) (1971), 41--44.%
\newline
\noindent \lbrack Hel] H. Helson, \textsl{Harmonic Analysis} (2$^{\text{nd}}$
ed.) Hindustan Book Agency, 1995. \newline
\noindent \lbrack HJ] J. Hoffmann-J\o rgensen, \textsl{Automatic Continuity,}
in: C. A. Rogers, J. Jayne, C. Dellacherie, F. Tops\o e, J. Hoffmann-J\o %
rgensen, D. A. Martin, A. S. Kechris, A. H. Stone, \textsl{Analytic sets.}
Academic Press, 1980; Part 3.2 .\newline
\noindent \lbrack Jab] E. Jab\l o\'{n}ska, On solutions of some
generalizations of the Go\l \c{a}b-Schinzel equation. \textsl{Functional
equations in mathematical analysis} (eds. T. M. Rassias and J. Brzd\k{e}k),
509--521, Springer, 2012.\newline
\noindent \lbrack KahM] P. Kahlig, J. Matkowski, A modified Go\l \k{a}%
b-Schinzel equation on a restricted domain (with applications to meteorology
and fluid mechanics), \textsl{\"{O}sterreich. Akad. Wiss. Math.-Natur. Kl.
Sitzungsber. II}, 211 (2002), 117-136 (2003).\newline
\noindent \lbrack Kor] J. Korevaar, \textsl{Tauberian theorems: A century of
development}. Grundl. math. Wiss. \textbf{329}, Springer, 2004.\newline
\noindent \lbrack Kuc] M. Kuczma, \textsl{An introduction to the theory of
functional equations and inequalities. Cauchy's equation and Jensen's
inequality.} 2nd ed., Birkh\"{a}user, 2009 [1st ed. PWN, Warszawa, 1985].%
\newline
\noindent \lbrack Moh] T. T. Moh, On a general Tauberian theorem. \textsl{%
Proc. Amer. Math. Soc.} \textbf{36} (1972), 167-172.\newline
\noindent \lbrack NgW] S. Ng, and S. Warner, Continuity of positive and
multiplicative functionals, \textsl{Duke Math. J.} \textbf{3}9 (1972),
281--284.\newline
\noindent \lbrack Ost1] A. J. Ostaszewski, Regular variation, topological
dynamics, and the Uniform Boundedness Theorem, \textsl{Topology Proceedings}%
, \textbf{36} (2010), 305-336.\newline
\noindent \lbrack Ost2] A.J. Ostaszewski, Beurling regular variation, Bloom
dichotomy, and the Go\l \k{a}b-Schinzel functional equation, \textsl{Aequat.
Math.}, Online First DOI 10.1007/s00010-014-0260-z.

\noindent \lbrack Oxt] J. C. Oxtoby: \textsl{Measure and category}, 2nd ed.
Graduate Texts in Math. \textbf{2}, Springer, 1980.\newline
\noindent \lbrack Pet] G. E. Petersen, Tauberian theorems for integrals II.
\textsl{J. London Math. Soc.} \textbf{5} (1972), 182-190. \newline
\noindent \lbrack Sen] E. Seneta, An interpretation of some aspects of
Karamata's theory of regular variation. \textsl{Publ. Inst. Math. (Beograd)}
(N.S.) \textbf{15 }(29) (1973), 111--119.\newline

\bigskip

\noindent Mathematics Department, Imperial College, London SW7 2AZ;
n.bingham@ic.ac.uk \newline
Mathematics Department, London School of Economics, Houghton Street, London
WC2A 2AE; A.J.Ostaszewski@lse.ac.uk\newpage

\end{document}